\documentclass[12pt]{article}
= 0.0in \topmargin = 0.3in \oddsidemargin=0.1in
\def\Dj{\hbox{D\kern-.73em\raise.30ex\hbox{-}
\raise-.30ex\hbox{}}}
\def\dj{\hbox{d\kern-.33em\raise.80ex\hbox{-}
\raise-.80ex\hbox{\kern-.40em}}}

\usepackage{amsmath}
\usepackage{mathrsfs}
\usepackage{amssymb}
\usepackage{amsfonts}
\usepackage{amsfonts,epsfig}
\usepackage{indentfirst}
\newtheorem{Theorem}{Theorem}[section]

\newtheorem{Remark}{Remark}[section]
\newtheorem{Lemma}{Lemma}[section]
\newtheorem{Corollary}{Corollary}[section]

\newenvironment {Proof} {\noindent {\bf Proof.}}{\hspace*{\fill}$\Box$\par\vspace{4mm}}

\newtheorem{theoremalph}{Theorem}

\allowdisplaybreaks[4]

\begin{document}

\title{The average distance of spanning trees in terms of independence number}

\author{Zhibin Du$^{a}$, Xuli Qi$^{b,}$\footnote{Corresponding author.}\\
\vspace{1mm}\\
$^a${\it School of Artificial Intelligence, South China Normal University}, \\
{\it Foshan 528225, P.R. China} \\
$^b${\it Department of Mathematics, Hebei University of Science and Technology}, \\
{\it Shijiazhuang 050018, P.R. China}
}

\date{}
\maketitle

\makeatletter

\def\@makefnmark{}

\makeatother

\footnotetext{   {\it E-mail addresses:} {\tt zhibindu@126.com}
(Z. Du), {\tt xuliqi@hebust.edu.cn} (X. Qi)}

\noindent {\bf Abstract }

Let $G$ be a connected graph with vertex set $V(G)$, and denote by $d_G(u,v)$ the distance from $u$ to $v$ in $G$, for any $u,v \in V(G)$.
The average distance of an $n$-vertex connected graph $G$, denoted by $\mu(G)$, is defined to be the average of all distances between all pairs of vertices in $G$, i.e.,
$\mu (G) = \binom{n}{2}^{-1} \sum_{\{u,v\} \subset V(G)}d_G(u,v)$.
The problem of
finding a spanning tree of minimum average distance is known to be
NP-hard, so establishing an upper bound for the minimum average distance among all spanning trees is of particular interest. Mukwembi (J. Graph Theory, 2014) showed that if $G$ is a connected graph of order $n$ with independence number $\alpha$, where $n > 2 \alpha - 1$, then $G$ has a spanning tree $T$ such that $\mu(T) \le \alpha + 2$.
In this paper, we first improve the upper bound to $\mu(T) < \alpha + 1$ for $\alpha \ge 1$, and then we find the bound could be further improved when $\alpha$ becomes larger, so a better upper bound
$$
\renewcommand{\arraystretch}{1.2}
\mu(T) < \left\{\begin{array}{ll}
\alpha + 1  & \mbox{if $1 \le \alpha \le 6$},\\
 \alpha + \frac{1}{2} + \frac{4(\alpha-1)}{\alpha^2}
 & \mbox{if $\alpha \ge 7$},
\end{array}
\right.
$$
is established later. In the end, we give a remark to indicate our new upper bound is best possible in the sense of asymptotics (when $n$ and $\alpha$ are large enough).
\\ \\


\section{Introduction}

Let $G$ be a connected graph with vertex set $V(G)$ throughout this paper. For $u,v \in V(G)$, denote by $d_G(u,v)$ the distance from $u$ to $v$ in $G$. The research related to distance of graphs has a long history in graph theory \cite{BuHa90}.

The average distance of an $n$-vertex connected graph $G$, denoted by $\mu(G)$, is defined to be the average of all distances between all pairs of vertices in $G$, i.e.,
$$
\mu (G) = \binom{n}{2}^{-1} \sum_{\{u,v\} \subset V(G)}d_G(u,v).
$$
Establishing the relationship between the average distance and other graph invariants has received a lot of attention in graph theory, and a large number of beautiful results can be found in many publications.

Among these, Dankelmann \cite{Dan94,Dan97} gave an upper and lower bounds for average distance of graphs in terms of independence number, matching number, and domination number, respectively, and he also characterized the graphs attaining these bounds. Kouider and Winkler \cite{KoWi97} showed that $\mu (G) \le \frac{n}{\delta + 1} + 2$ for all connected graphs $G$ of order $n$ with minimum degree $\delta$.
Dankelmann {\it et al.} \cite{DaMuSw09}
investigated the relationship between the average distance and vertex-connectivity, and proved that $\mu (G) \le \frac{n}{2} (\kappa + 1) + O(1)$, when $G$ is a $\kappa$-vertex-connected graph of order $n$, where $\kappa \ge 3$ is odd.
Hansen {\it et al.} \cite{HaHeKiMaSc09} showed that $\mu (G) \le \frac{1}{2} F(G)$, where $F(G)$ denotes the maximum order of an induced forest of $G$.
Recently, Mukwembi \cite{Muk24} gave an upper bound on the average distance of a connected graph of given order and minimum degree where irregularity index is prescribed, strengthening the results in \cite{DaEn00,KoWi97}.

It is worth mentioning that the definition of average distance $\mu (G)$ is exactly equivalent to that of the well-known distance-based graph invariant--Wiener index $W(G)$, when the order of $G$ is fixed. The Wiener index of a connected graph $G$ is defined as the sum of all distances between all pairs of vertices in $G$, see \cite{Win47}, i.e.,
$$
W (G) = \sum_{\{u,v\} \subset V(G)}d_G(u,v).
$$
Clearly,
$$
W(G) = \binom{n}{2} \mu(G),
$$
where $n$ is the order of $G$.
As a consequence, the results and applications of average distance and Wiener index are equivalent.

The Wiener index (also called Wiener number) is commonly regarded as the first proposed topological index in chemical graph theory \cite{GuPo86}, which is an effective numerical descriptor to measure a variety of properties of molecular compounds, especially, boiling points, heats of formation, heats of vaporization,
molar volumes, and molar refractions \cite{NiTrMi95}.

Besides the remarkable applications in chemistry, it also found some practical applications, such as communication, facility location, etc. In particular, a smaller Wiener index would lead to a lower cost for communications and facility locations. So for the sake of saving cost, finding a spanning tree with minimum Wiener index (average distance) for a given graph is of particular interest. However, unfortunately, this problem turned out to be NP-hard \cite{JoLeKa78}. Since it is not easy to search optimal spanning tree (with minimum average distance), we may try to find its best bound (if exist) instead. Therefore it is meaningful to find some graph invariants to bound the minimum average distance among spanning trees.

For more results about average distance (Wiener index) of graphs, we refer the readers to papers \cite{BeRiSm01,Chen23,Da12,DaOeWe04,DePeWa09,DuZh10,Hen86,YeGu94}, and two surveys \cite{DoEnGu01,Ple84}.

A pioneering research for finding a spanning tree with minimum Wiener index (average distance), under the background of total cost of communication, was done by Hu in 1974, see \cite{Hu74}. Entringer {\it et al.} \cite{EnKlSz96} found a spanning tree $T$ of a connected graph $G$, such that $\mu(T) \le 2\mu (G)$. Subsequently,
Dankelmann and Entringer \cite{DaEn00} proved that every connected graph of order $n$ with minimum degree $\delta$ has a spanning tree with average distance at most $\frac{n}{\delta + 1} + 5$.
Ten years ago, Mukwembi \cite{Muk14} reported an upper bound of minimum average distance of spanning trees in terms of independence number.

Motivated by \cite{Muk14}, this paper will further contribute to the research about the relationship between average distance of spanning trees and independence number. An independent set is a vertex subset in which no pair is adjacent.
The independence number of a graph $G$, denoted by $\alpha(G)$, or simply $\alpha$, is
the maximum cardinality of an independent set of $G$.

Initially, this topic originated from an interesting conjecture, posed by a computer program called Graffiti, which states $\mu (G) \le \alpha (G)$ holds for all connected graphs $G$, see \cite{FaWa86,FaWa87}. In the same paper \cite{FaWa86}, a somewhat weaker version $\mu (G) < \alpha (G) + 1$ was solved. Two years later, Chung \cite{Chu88} completed the solution of this conjecture. Later, Dankelmann \cite{Dan94} extended this result to determine the unique graph with maximum average distance with fixed independence number.
Ten years ago, Mukwembi \cite{Muk14} used the independence number to bound the minimum average distance among all the spanning trees.
More precisely, his main result is listed as follows.

\begin{theoremalph} \cite{Muk14} \label{old}
Let $G$ be a connected graph of order $n$ and independence number $\alpha$. Then $G$ has a spanning tree $T$ with
$$
\renewcommand{\arraystretch}{1.2}
\mu(T) \le \left\{\begin{array}{ll}
 \alpha + 2
 & \mbox{if $n >2 \alpha - 1$},\\
\frac{2}{3} \alpha
 & \mbox{if $n \le 2 \alpha - 1$.}
\end{array}
\right.
$$
\end{theoremalph}

The upper bound $\frac{2}{3} \alpha$ when $n \le 2 \alpha - 1$ is attained by the path trivially, thus such case is likely to be of less interest to us.
In this paper, our object is the other case $n >2 \alpha - 1$, we will improve the inequality $\mu(T) \le \alpha + 2$, described in Theorem \ref{old} when $n >2 \alpha - 1$, to $\mu(T) < \alpha + 1$ when $\alpha \ge 1$, and to
$$
\mu(T) < \alpha + \frac{1}{2} + \frac{4(\alpha-1)}{\alpha^2}
$$
when $\alpha \ge 7$. Note that $\frac{4(\alpha-1)}{\alpha^2} \rightarrow 0$ as $\alpha \rightarrow \infty$, so the extremity of this bound is $\alpha + \frac{1}{2}$. It is worth mentioning that the extreme upper bound $\alpha + \frac{1}{2}$ is best possible in the sense of asymptotics (when $n$ and $\alpha$ are large enough).

Our paper is organized as follows: In Sections $2$ and $3$, some useful definitions and lemmas are introduced, then the improvements of the upper bound, from $\alpha + 2$ to $\alpha + 1$ when $\alpha \ge 1$ and further to
$$
\renewcommand{\arraystretch}{1.2}
\mu(T) < \left\{\begin{array}{ll}
\alpha + 1  & \mbox{if $1 \le \alpha \le 6$},\\
 \alpha + \frac{1}{2} + \frac{4(\alpha-1)}{\alpha^2}
 & \mbox{if $\alpha \ge 7$},
\end{array}
\right.
$$
would be established in the following two sections, respectively, finally we  claim that the improved upper bound  is  best possible in the sense of asymptotics (when $n$ and $\alpha$ are large enough).

\section{Preliminaries}

A vertex of degree $1$ is said to be a pendent vertex.

\begin{figure}[h]
  \center
  \includegraphics [width = 12cm]{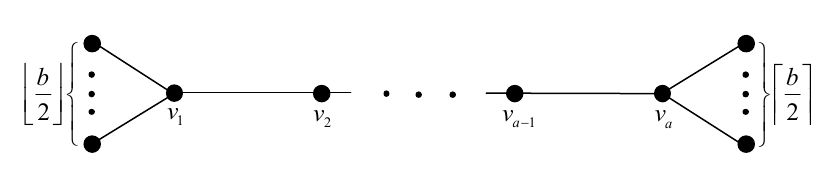}\vspace{0.2cm}
  \caption {The tree $H_{a,b}$. } \label{fig-key}
\end{figure}

Let $H_{a,b}$ be the tree of order $a+b$ obtained from the path $P_a = v_1 v_2 \dots v_a$ by attaching $\lfloor\frac{b}{2} \rfloor$ pendent vertices to $v_1$, and $\lceil\frac{b}{2} \rceil$ pendent vertices to $v_a$, where $a \ge 2$ and $b \ge 2$, see Figure \ref{fig-key}.

The tree $H_{a,b}$ plays an important role in the researches about extremal values of many distance-based graph invariants, e.g., Wiener index \cite{Dan94,Du12,Ent99,Shi93}, distance spectral radius \cite{DuIiFe13}.

For an edge subset $M$ of the graph $G$, $G-M$ denotes the
graph obtained from $G$ by deleting the edges in $M$, and for an edge subset $M^*$ of the complement of $G$, $G+M^*$
denotes the graph obtained from $G$ by adding edges in $M^*$. For simplicity, if $M$ contains a single edge, say $M = \{uv\}$, we use $G- uv$ instead of $G- \{uv\}$ for short. Similarly, if $M^*$ contains a single edge, say $M^* = \{uv\}$, we use $G + uv$ instead of $G + \{uv\}$ for short. Moreover, for $v \in V(G)$, let $G-v$ be the graph obtained from $G$ by deleting the vertex $v$ and its incident edges.

\section{Lemmas}

In this section, we will present some useful lemmas for our main results.

First of all, we consider the effects on Wiener indices under two graph transformations.
Let $N_G(v)$ denote the set of neighbors of $v$ in $G$.

\begin{figure}[h]
  \center
  \includegraphics [width = 10cm]{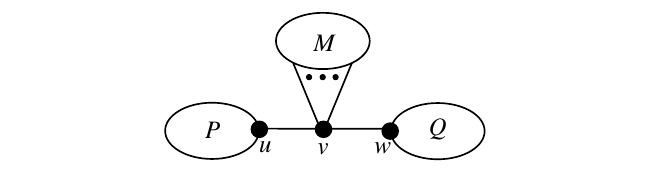}\vspace{0.2cm}
  \caption {The tree $T$ in Lemma \ref{le-tran}. } \label{fig-tran}
\end{figure}

\begin{Lemma} \label{le-tran}
Let $T$ be a tree of the form as depicted in Figure \ref{fig-tran}. Assume that $v_1, v_2, \dots, v_t$ are the neighbors of $v$ lying in $M$, where $t \ge 1$.
Consider
$$
T_u = T- \{v v_1, v v_2, \dots, v v_t\} + \{u v_1, u
v_2, \dots, u v_t\}
$$
and
$$
T_w = T- \{v v_1, v v_2, \dots, v v_t\} + \{w v_1, w v_2, \dots, w
v_t\}.
$$
Then
$$
W(T) <
\max\{W(T_u), W(T_w)\}.
$$
\end{Lemma}

\begin{Proof}
If $|V(P)| \ge |V(Q)|$, then
$$
W(T_w) - W(T) = |V(M)| (|V(P)| - |V(Q)| + 1) > 0.
$$
If $|V(P)| < |V(Q)|$, then similarly
$$
W(T_u) - W(T) = |V(M)| (|V(Q)| - |V(P)| + 1) > 0.
$$
The result holds clearly.
\end{Proof}

\begin{figure}[h]
  \center
  \includegraphics [width = 10cm]{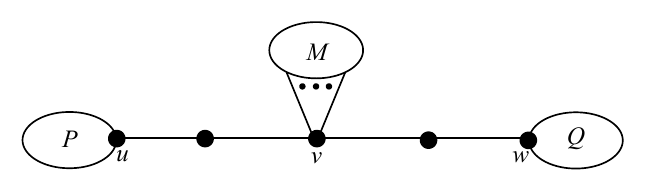}\vspace{0.2cm}
  \caption {The tree $T$ in Remark \ref{remark-add}. } \label{fig-tran-1}
\end{figure}

\begin{Remark} \label{remark-add}
If we insert one vertex between $v$ and $u$, and another vertex between $v$ and $w$, the tree under consideration is as depicted in Figure \ref{fig-tran-1}, then the conclusion $W(T) <
\max\{W(T_u), W(T_w)\}$ is still valid, following an analogous proof as above.
\end{Remark}

Let us recall the characterization of the unique tree with maximum Wiener index among trees with given order and number of pendent vertices.

\begin{Lemma} \cite{Ent99,Shi93} \label{le-pendent}
Let $T$ be a tree of order $n \ge 4$ with $p$ pendent vertices, where $2 \le p \le n - 2$. Then
$$
W(T) \le \frac{1}{12} (2n^3 - (3p^2-6p+2)n+p^3+3p^2-10p),
$$
with equality if and only if $T \cong H_{n-p,p}$ with even $p$. 
\end{Lemma}

\begin{Remark} \label{le-pendent-com}
It is easy to verify that the upper bound
$$
\frac{1}{12} (2n^3 - (3p^2-6p+2)n+p^3+3p^2-10p)
$$
on $W(T)$ claimed in Lemma \ref{le-pendent} is decreasing for $p \ge 1$ (by checking its derivative with respect to $p$ is negative).
\end{Remark}

Let $\mathcal{T}_{a,b}$ be the set of trees of order $a+b$ obtained from a tree $H$ of order $a$ by attaching $b$ pendent vertices to some vertices of $H$, where $a \ge 2$ and $b \ge 1$.

Mukwembi \cite{Muk14} established a tight upper bound for the Wiener indices of trees in $\mathcal{T}_{a,b}$, which plays a crucial role in proving his main result (Theorem \ref{old}).
We found that it can also be proved in an alternative and shorter way, by using the number of pendent vertices.

\begin{Lemma} \cite{Muk14} \label{old-le}
Let $T \in \mathcal{T}_{a,b}$, where $a \ge 2$ and $b \ge 1$. Then
$$
W (T) \le  \frac{1}{12}(2 a ^3 + 6a^2 b + 3a b^2 + 6a b -2 a + 9 b^2 - 12 b  ),
$$
with equality if and only if $T \cong H_{a,b}$ with even $b$.
\end{Lemma}

\begin{Proof}
Let $T$ be a tree in $\mathcal{T}_{a,b}$. Denote by $p$ the number of pendent vertices in $T$. From the definition of trees in $\mathcal{T}_{a,b}$, obviously $p \ge b$.

Let $n=a+b$, and set
$$
f(n,x) := \frac{1}{12} (2n^3 -(3x^2-6x+2)n+x^3+3x^2-10x),
$$
where $1 \le x \le n- 2$.

If $p = b$, then by Lemma \ref{le-pendent}, we have
$$
W(T) \le f(n,p) = f(n,b),
$$
with equality if and only if $T \cong H_{a,b}$ with even $b$.

If $p > b$, then by Lemma \ref{le-pendent} and Remark \ref{le-pendent-com}, it holds that
$$
W(T) \le f(n,p) < f(n,b).
$$

Note that
$$
f(n,b) = f(a+b,b) = \frac{1}{12}(2 a ^3 + 6a^2 b + 3a b^2 + 6a b -2 a + 9 b^2 - 12 b  ).
$$
Now the result follows.
\end{Proof}

For convenience, sometimes we would use an equivalent form of Lemma \ref{old-le} as follows.

\begin{Lemma}  \label{old-le-eq}
Let $T \in \mathcal{T}_{t,n-t}$, where $2 \le t \le n - 1$. Then
$$
W (T) \le  \frac{1}{12}( -t^3 + 3t^2 + (3n^2 - 12n + 10)t + 9n^2 - 12n),
$$
with equality if and only if $T \cong H_{t,n-t}$ with even $n-t$.
\end{Lemma}

An equivalent form of Remark \ref{le-pendent-com} is also presented here.

\begin{Remark} \label{remark}
The upper bound
$$
\frac{1}{12}( -t^3 + 3t^2 + (3n^2 - 12n + 10)t + 9n^2 - 12n)
$$
about $W(T)$ claimed in Lemma \ref{old-le-eq} is increasing for $t \ge 2$ (by checking its derivative with respect to $t$ is positive).
\end{Remark}

We further expound how the set $\mathcal{T}_{a,b}$ contributes to the proof of Theorem \ref{old} as in \cite{Muk14}.
Let $G$ be a connected graph of order $n$ and independence number $\alpha$. In the proof of Theorem \ref{old} in \cite{Muk14}, it was shown that there exists a subtree $\bar{T}$ of $G$ of order $t$ such that the remaining $n-t$ vertices are adjacent to some vertices in $\bar{T}$, where $t \le 2 \alpha - 1$. That is to say, we can find a spanning tree $T$ of $G$ such that $T \in \mathcal{T}_{t, n-t}$ with $t \le 2 \alpha - 1$.
Here we would like to strengthen it in the aspect of the lengths of pendent paths and  internal paths for our purpose.

Let $d_G(v)$ denote the degree of the vertex $v$ in $G$. Suppose that $P= v_0 v_1 \dots v_k$ is a path in $G$ of length $k \ge 1$, where $d_G(v_0) \ge 3$, $d_G(v_1) = d_G(v_2) = \dots = d_G(v_{k-1}) = 2$.
If $d_G(v_k) \ge 3$, then $P$
is said to be an internal path of $G$. If $d_G(v_k) = 1$, then $P$ is said to be a pendent path of $G$.

\begin{Lemma} \label{le-final-add}
Let $G$ be a connected graph of order $n$ and independence number $\alpha$, where $n \ge 2 \alpha$. Then we can find a spanning tree $T$ of $G$ such that $T \in \mathcal{T}_{t, n-t}$, for $t \le 2 \alpha - 1$. In particular, when $t = 2 \alpha - 1$, say $T$ is formed from $\bar{T}$ of order $2 \alpha - 1$ by attaching $n - 2 \alpha + 1$ pendent vertices to some vertices of $\bar{T}$,
\begin{itemize}

\item[(i)]
when $\bar{T}$ is not a path, every pendent path or internal path of $\bar{T}$ is of even length;

\item[(ii)]

when $\bar{T}$ is a path, say $\bar{T} = v_1 v_2 \dots v_{2 \alpha - 1}$, the $n - 2 \alpha + 1$ pendent vertices can only be attached to the vertices among $v_1, v_3, \dots, v_{2 \alpha - 1}$.

\end{itemize}
\end{Lemma}

\begin{Proof}
First of all, let us recall the arguments about the existence of spanning tree $T$ of $G$ such that $T \in \mathcal{T}_{t, n-t}$, for some $t \le 2 \alpha - 1$, which was given in \cite{Muk14}. The desired spanning tree $T$ can be inductively constructed as follows:
\begin{itemize}

\item
First, choose an arbitrary vertex $v$ of $G$,  let $A_1 = \{v\}$, and
$T_1$ be the tree with only one vertex $v$, it is  a subgraph of $T$.

\item

Now assume that the set $A_{i-1}$ and the tree $T_{i-1}$ have been constructed, where $i \ge 2$.
If there exists a vertex $x \not\in V(T_{i-1})$ such that $x$ is not adjacent to any vertex in $A_{i-1}$, and $d_G(x, y) = 2$ for some $y \in A_{i-1}$, let $A_i = A_{i-1} \cup \{x\}$. Moreover, let
$xzy$ be a shortest path between $x$ and $y$ in $G$, where $z\in V(T_{i-1})$ is allowed, we define $T_i = xzy \cup T_{i-1}$.

\item
Repeat this process until it is exhausted, at that time we get the set $A_k$ and the tree $T_k$, every vertex $w \not\in V(T_k)$ would be adjacent to some vertex, write $w'$, in $A_k$.
Meanwhile, a spanning tree $T$ of $G$ with edge set $E(T_k) \cup \{ww' : w \not\in V(T_k), w' \in A_k\}$ can be found.
That is to say, $T \in \mathcal{T}_{t, n-t}$, where $t$ is the order of $T_k$.
\end{itemize}

On one hand, it is observed that one or two vertices would be added when we form $T_i$ from $T_{i-1}$, for any $2 \le i \le k$, which implies that $t \le 1 + 2 (k - 1) = 2k - 1$. On the other hand, the vertices of $A_k$ are pairwise independent, thus $A_k$ is an independent set of $G$, leading to $k \le \alpha$. Combining the two inequalities, it results in $t \le 2\alpha - 1$.

In particular, when $t = 2\alpha - 1$, we should require that $k = \alpha$, and each time exactly two vertices are added to form $T_i$ from $T_{i-1}$, for any $2 \le i \le k$. More importantly, the latter requirement infers the following two observations:
\begin{itemize}

\item

When $T_k$ is the path on $2\alpha - 1$ vertices, say $T = v_1 v_2 \dots v_{2 \alpha - 1}$, we have $A_{k} = \{v_1, v_3, \dots, v_{2 \alpha - 1}\}$, thus the spanning tree $T$ of $G$ is obtained from $T_k$ by attaching some pendent vertices  to the vertices in  $A_{k}$.

\item

When $T_k$ is not a path, there exist some vertices of degrees at least $3$ in $T_k$, all such vertices must belong to $A_k$, thus every pendent path or internal path of $T_k$ should be of even length.

\end{itemize}

The desired result follows easily.
\end{Proof}

\section{Preliminary refinement of Theorem \ref{old}}

In this section, we will claim that the upper bound $\alpha + 2$ presented in \cite[Theorem 2.1]{Muk14} (i.e., Theorem \ref{old}), when $n > 2 \alpha - 1$, can be improved to $\alpha + 1$.

As it is stated in Lemma \ref{le-final-add}, when $n \ge 2 \alpha$, we may construct a spanning tree $T$ such that $T \in \mathcal{T}_{t, n-t}$, for $t \le 2 \alpha - 1$.
Combining with Lemma \ref{old-le-eq} and Remark \ref{remark}, we can get

\begin{eqnarray*}
W(T) &\le& \frac{1}{12}( -t^3 + 3t^2 + (3n^2 - 12n + 10)t + 9n^2 - 12n) \\
&\le& \frac{1}{12}( -(2 \alpha - 1)^3 + 3(2 \alpha - 1)^2 + (3n^2 - 12n + 10)(2 \alpha - 1) + 9n^2 - 12n) \\
&=& \frac{1}{6} (-4 \alpha^3 + 12 \alpha^2 + (3 n^2 -12 n + 1) \alpha +3 n^2 - 3).
\end{eqnarray*}
In particular, all the equalities hold
if and only if $T \cong H_{2 \alpha - 1, n - 2 \alpha + 1}$ with even $n - 2 \alpha + 1$.
Then we get
$$
\mu (T) - \alpha \le \frac{-4 \alpha^3 + 12\alpha^2 - (9n - 1)\alpha + 3n^2 - 3}{3 n (n-1)}.
$$
It is easily verified that
$$
\frac{-4 \alpha^3 + 12\alpha^2 - (9n - 1)\alpha  + 3n^2 - 3}{3 n (n-1)}
$$
is decreasing for $\alpha \ge 1$ (by checking its derivative with respect to $\alpha$ is negative), thus
$$
\mu (T) - \alpha \le \frac{-4 + 12 - (9n - 1) + 3n^2 - 3}{3 n (n-1)} = \frac{3n^2 - 9 n + 6}{3 n (n-1)} < 1,
$$
i.e.,
$$
\mu (T) < \alpha + 1.
$$

Now it leads to our first improvement of Theorem \ref{old}.

\begin{Theorem} \label{th-new}
Let $G$ be a connected graph of order $n$ and independence number $\alpha$, where $n \ge 2 \alpha$ and $\alpha \ge 1$. Then $G$ has a spanning tree $T$ with
$$
\mu(T) < \alpha + 1.
$$
\end{Theorem}

Clearly, $\mu(G) \le \mu(T)$, for any spanning tree $T$ of $G$. The following corollary can be obtained immediately.

\begin{Corollary} \label{con}
Let $G$ be a connected graph of order $n$ and independence number $\alpha$, where $n \ge 2 \alpha$ and $\alpha \ge 1$. Then
$$
\mu(G) < \alpha + 1.
$$
\end{Corollary}

\section{Further improvement for the upper bound}

It is obvious that the complete graph $K_n$ is the unique $n$-vertex connected graph with independence number $\alpha = 1$, and its spanning tree with minimum average distance is exactly the star $S_n$. Note that
$$
\mu (S_n) = 2 \left( 1 - \frac{1}{n} \right),
$$
which approaches to $\alpha + 1 = 2$ for large $n$. It means that the upper bound $\alpha + 1$ is best possible for $\alpha = 1$, when $n$ is sufficiently large.

Now it naturally raises a further problem: Is it possible to improve the upper bound $\alpha + 1$ when $\alpha$ is somewhat larger?

In this section, we will go ahead along this direction, and manage to improve the upper bound to
$$
\renewcommand{\arraystretch}{1.2}
\mu(T) < \left\{\begin{array}{ll}
\alpha + 1  & \mbox{if $1 \le \alpha \le 6$},\\
 \alpha + \frac{1}{2} + \frac{4(\alpha-1)}{\alpha^2}
 & \mbox{if $\alpha \ge 7$},
\end{array}
\right.
$$
by exploring more detailed structures of the graph and its spanning trees.

In view of Lemma \ref{le-final-add}, our interest should be focused on the trees in $\mathcal{T}_{t,n-t}$ with $t \le 2 k- 1$.
Our first result deals with the trees in $\mathcal{T}_{t,n-t}$, for $t \le 2 k- 2$.

\begin{Lemma} \label{le-add}
Let $T \in \mathcal{T}_{t,n-t}$, where $t \le 2 k- 2$, $2 \le t \le n - 2$, and $k \ge 2$.
Then
$$
\mu (T) < k + \frac{1}{2}.
$$
\end{Lemma}

\begin{Proof}
For fixed $n$, let
$$
f(t) = \frac{1}{12}( -t^3 + 3t^2 + (3n^2 - 12n + 10)t + 9n^2 - 12n),
$$
where $t \ge 2$.
From Lemma \ref{old-le-eq}, we know that
$W(T) \le f(t)$.
And $f(t)$ is increasing for $t \ge 2$, from Remark \ref{remark}.

Since $t \le 2 k - 2$, now we have
$$
W(T) \le f(t) \le f(2 k - 2).
$$
Notice that
$$
f(2k -2) = \frac{1}{12} (-8k ^3 + 36 k^2 + 2( 3n^2 -12 n - 14 ) k  + 3 n^2 + 12n  ).
$$

Furthermore, we get
\begin{eqnarray*}
&&\mu(T) - k \nonumber \\
&\le& \binom{n}{2}^{-1} f(2 k - 2) - k \nonumber \\
&=&  \frac{-8k ^3 + 36 k^2 - 2( 9 n + 14 ) k  + 3 n^2 + 12n }{6 n (n-1)}.
\end{eqnarray*}
It is easily verified that
$$
\frac{-8k ^3 + 36 k^2 - 2( 9 n + 14 ) k  + 3 n^2 + 12n }{6 n (n-1)}
$$
is decreasing for $k \ge 2$. So
\begin{eqnarray*}
\mu(T) - k &\le& \frac{-8\cdot 2 ^3 + 36 \cdot 2^2 - 2( 9 n + 14 ) 2  + 3 n^2 + 12n }{6 n (n-1)} \\
&=& \frac{  n^2 - 8 n + 8}{2 n (n-1)} \\
&<& \frac{1}{2},
\end{eqnarray*}
i.e.,
$$
\mu (T) < k + \frac{1}{2}.
$$

The result follows.
\end{Proof}

In the following several lemmas, we will investigate the trees in $\mathcal{T}_{2k-1,n-2k+1}$ with some particular structures.

\begin{Lemma} \label{diameter}
Let $T$ be a tree obtained from the path $P_{2k-1}$ by attaching $b \ge 1$ pendent vertices to some vertices of $P_{2k-1}$, where $k \ge 2$. If the diameter of $T$ is $2k-2$ or $2k-1$, then
$$
\mu (T) < k + \frac{1}{2}.
$$
\end{Lemma}

\begin{Proof}
Clearly, the number of pendent vertices in $T$ is $b+1$ or $b+2$ (notice that if the number of pendent vertices is $b$, then the diameter of $T$ is $2k$). By Lemma \ref{le-pendent} and Remark \ref{le-pendent-com}, we know that the maximum Wiener index is attained by $H_{2k-2,b+1}$ (the number of pendent vertices is $b+1$), more precisely, by setting $n = 2k -1+b$ (so $b + 1 = n - 2k + 2$), we have
\begin{eqnarray*}
W(T) &\le& \frac{1}{12} (2n^3 - (3(n - 2k + 2)^2-6(n - 2k + 2)+2)n  \\
&&+(n - 2k + 2)^3+3(n - 2k + 2)^2-10(n - 2k + 2))  \\
&=& \frac{1}{12} (-8k ^3 + 36 k^2 + 2( 3n^2 -12 n - 14 ) k  + 3 n^2 + 12n  ).
\end{eqnarray*}
In the proof of Lemma \ref{le-add}, we have shown that
$$
\mu (T) < k + \frac{1}{2},
$$
as desired.
\end{Proof}

\begin{Remark} \label{remark-final-add}
In fact, the above proof indicates that when the number of pendent vertices of $T$ is $n - 2k + 2$,
$$
\mu (T) < k + \frac{1}{2}.
$$
\end{Remark}

\begin{Lemma} \label{le-final-2}
Let $T$ be a tree obtained from the path $P_{2k} = v_1 v_2 \dots v_{2k}$ by attaching $n-2k$ pendent vertices to some vertices of $P_{2k}$, except $v_2$, where $k \ge 2$. If there are at most $\frac{n-2k+1}{k}$ pendent vertices attached to $v_1$, then
$$
\mu (T) < k + \frac{1}{2} + \frac{4(k-1)}{k^2}.
$$
\end{Lemma}

\begin{Proof}
In virtue of Lemma \ref{le-tran}, we may assume that all the pendent vertices outside $P_{2k}$ can only be attached to $v_1, v_3, v_{2k}$.
Denote by $a$ the number of pendent vertices attached to $v_1$, and $b$ the number of pendent vertices attached to $v_{2k}$, where $a \le \frac{n-2k+1}{k}$, see Figure \ref{fig-final-add}.
Then we have
\begin{eqnarray} \label{le-w-expression-next}
W (T) &=& (n-1)^2 + a(n - a - 2) + (a + 1)(n - a - 3) \nonumber \\
&&+ \sum_{i = 0}^{2k-4}(b + i)(n - b - i - 2) \nonumber\\
&=& -2 a^2 + 2(n-3)a - (2k - 3) b^2 + (2k - 3) (n - 2k + 2) b \nonumber\\
&& + n^2 +(2k-5)(k-1)n - \frac{1}{3} k (8k^2 -30k + 31).
\end{eqnarray}

\begin{figure}[h]
  \center
  \includegraphics [width = 12cm]{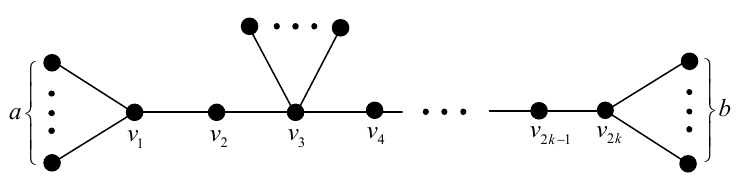}\vspace{0.2cm}
  \caption {The tree $T$ used in the proof of Lemma \ref{le-final-2}. } \label{fig-final-add}
\end{figure}

It is not hard to verify that $W(T)$ is increasing in $a \le \frac{n-2k+1}{k}$ (since $\frac{n-2k+1}{k} \le \frac{n-3}{2}$), thus the maximum value of $W(T)$ is attained at $a = \frac{n-2k+1}{k}$. Substitute $a = \frac{n-2k+1}{k}$ into \eqref{le-w-expression-next}, the upper bound of $W(T)$ would be expressed in terms of a quadratic function on $b$ as
$W(T) \le \frac{f(b)}{3k^2}$,
where
\begin{eqnarray*}
f(b) &=& -3k^2(2k-3) b^2  + 3k^2(2k-3) (n - 2k + 2) b + 3(k^2 + 2k - 2)n^2 \\
&&+ 3 (2k^4 - 7 k^3 +k^2 +4k -4) n - 8 k^5 - 6 (2b - 5) k^4 - (6b^2-30b+31) k^3 \\
&&+ 3 (3b^2 - 6b + 4)k^2 + 6k - 6.
\end{eqnarray*}
The maximum value of $f(b)$ is attained at $b = \frac{n-2k+2}{2}$. That is to say, after substitute $a = \frac{n-2k+1}{k}$ and $b = \frac{n-2k+2}{2}$ into \eqref{le-w-expression-next}, we will obtain the maximum value of $W(T)$.

Further,
after a series of standard (but somewhat tedious) calculations, we are able to conclude with
\begin{eqnarray*}
&&\mu (T) - k \\
&=& \frac{W(T)}{\binom{n}{2}} - k \\
&\le& \frac{3 (k^2 + 8k - 8) n^2 - 6 (3k^3 + 4 k^2 - 8k + 8) n - 4 (2 k^5 - 9k ^4 + 7k^3 - 3 k^2 - 6k + 6)}{6k^2 n (n-1)}\\
&<& \frac{1}{2} + \frac{4(k-1)}{k^2},
\end{eqnarray*}
i.e.,
$$
\mu (T) < k + \frac{1}{2} + \frac{4(k-1)}{k^2}.
$$

The proof is completed.
\end{Proof}

\begin{Lemma} \label{le-extra}
Assume that $T$ is obtained from a tree $H_T$ of order $2k-1$ by attaching $n-2k+1$ pendent vertices to some vertices of $H_T$, i.e., $T \in \mathcal{T}_{2k-1,n-2k+1}$, where $n \ge 2k$ and $k \ge 4$. If $H_T$ is not a path, and each pendent path or internal path in $H_T$ is of even length,  then
$$
\mu (T) < k + \frac{1}{2} + \frac{1}{2(2k-5)}.
$$
\end{Lemma}

\begin{Proof}
Since our object is an upper bound of $\mu(T)$ (or equivalently $W(T)$), we may assume that $T$ is a tree with maximum $W(T)$ among the trees in $\mathcal{T}_{2k-1,n-2k+1}$ formed from a non-path $H_T$, and each pendent path or internal path in $H_T$ is of even length.
We claim that $H_T$ is the tree obtained by attaching two paths on two vertices to an end vertex of the path $P_{2k - 5}$.
For proving this claim, let us say $v_1 v_2 \dots v_t$ a diametrical path of $H_T$.

Suppose that there are two different indices $3 \le i,j \le t - 2$ such that $d_{H_T} (v_i), d_{H_T} (v_j) \ge 3$. By Remark \ref{remark-add}, we can construct another tree $T_1$ with $W(T_1) > W(T)$ by  switching all the neighbors of $v_i$ in $T$ different from $v_{i-1},v_{i+1}$ to be the neighbors of $v_{i-2}$ or $v_{i+2}$ in $T_1$. It is observed that $H_{T_1}$ is not a path either (because $d_{H_{T_1}} (v_j) \ge 3$), and each pendent path or internal path in $H_{T_1}$ is still of even length, leading to a contradiction to the maximality of $W(T)$.

We have known that there is exactly one index $3 \le i \le t - 2$ with $d_{H_T} (v_i) \ge 3$. Following the same reasoning as last paragraph, we can use Remark \ref{remark-add} again to deduce that $i = 3$ or $t - 2$, and $d_{H_T} (v_i) = 3$, that is to say,  $H_T$ is obtained by attaching two paths on two vertices to an end vertex of the path $P_{2k - 5}$.

Finally, we estimate the upper bound of $\mu(T)$, where $T$ is obtained from $H_T$ by attaching $n-2k+1$ pendent vertices to some vertices of $H_T$.
From Lemma \ref{le-tran}, we can assume that the $n-2k+1$ pendent vertices are attached to  the three end vertices of $H_T$. More precisely, we may assume that $T$ is of the form as depicted in Figure \ref{fig-add}, where $a + b + c = n - 2k + 1$, equivalently $a = n - 2k + 1 - b - c$.

\begin{figure}[h]
  \center
  \includegraphics [width = 9cm]{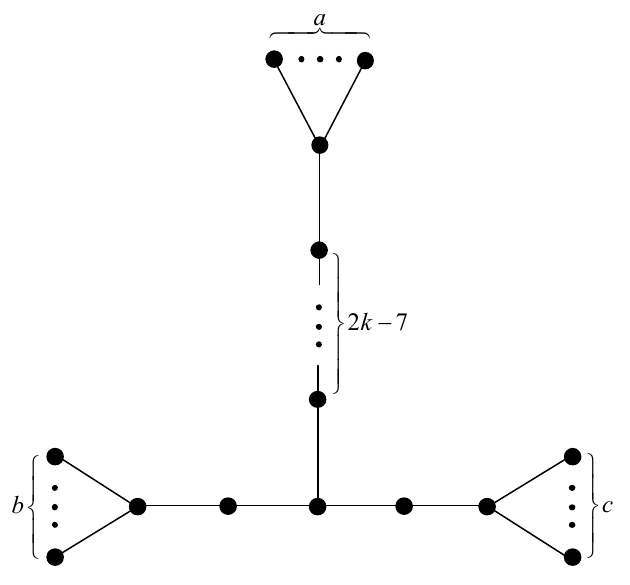}\vspace{0.2cm}
  \caption {The tree $T$ used in the proof of Lemma \ref{le-extra}. } \label{fig-add}
\end{figure}

The value of $W(T)$ is expressed as
\begin{eqnarray} \label{le-w-expression-new}
W (T) &=& (n-1)^2 + \sum_{i = 0}^{2k - 7}(a + i)(n - a - i - 2)+ b(n -b  - 2) + (b + 1)(n - b - 3) \nonumber\\
&&+ c(n - c  - 2) + (c + 1)(n - c - 3) \nonumber \\
&=& (n-1)^2 + \sum_{i = 0}^{2k -7}((n - 2k + 1 - b - c) + i)(n - (n - 2k + 1 - b - c) - i - 2)\nonumber\\
&&+ b(n -b  - 2) + (b + 1)(n - b - 3) + c(n - c  - 2) + (c + 1)(n - c - 3).
\end{eqnarray}
Analogous to the handling of \eqref{le-w-expression-next} in the proof of Lemma \ref{le-final-2} (in a more complicated form), the maximum value of $W(T)$ would be attained at
$$
b = c = \frac{(k-2)n - 2k^2 + 3k + 6}{2(2k -5 )}.
$$
After substitute such $b$ and $c$ into \eqref{le-w-expression-new}, and implement a series of standard (but somewhat tedious) calculations, we will finally get
$$
\mu (T) - k =  \frac{W(T)}{\binom{n}{2}} - k < \frac{1}{2} + \frac{1}{2(2k-5)},
$$
or equivalently,
$$
\mu (T) < k + \frac{1}{2} + \frac{1}{2(2k-5)}.
$$

The proof is completed.
\end{Proof}

\begin{Lemma} \label{le-key}
Let $G$ be a connected graph of order $n$ and independence number $\alpha$, where $n \ge 2 \alpha$ and $\alpha \ge 2$. Assume that $T$ is a spanning tree of $G$ obtained from the path $P_{2\alpha-1} = v_1 v_2 \dots v_{2\alpha-1}$ by attaching $n - 2\alpha + 1$ pendent vertices to some vertices among $v_1, v_3, \dots, v_{2\alpha-1}$.
Then there exists a spanning tree $T^*$ of $G$ such that
$$
\mu(T^*) < \alpha + \frac{1}{2} + \frac{4(\alpha-1)}{\alpha^2}.
$$
\end{Lemma}

\begin{Proof}
Clearly, $T \in \mathcal{T}_{2\alpha-1,n - 2\alpha + 1}$.
On the other hand, note that the diameter of $T$ is possibly $2\alpha - 2$, $2\alpha - 1$, or $2\alpha$. If its diameter is $2\alpha - 2$ or $2\alpha - 1$, then
from Lemma \ref{diameter},
$$
\mu (T) < \alpha + \frac{1}{2}.
$$
So we may assume that the diameter of $T$ is $2\alpha$, it implies that there are some pendent vertices attached to $v_1$ and $v_{2\alpha-1}$ in $T$. Clearly, $T$ has $n - 2\alpha + 1$ pendent vertices.

\noindent {\bf Case 1.} There is some edge joining $v_i$ and $v_j$ in $G$, where $|j - i| \ge 2$.

From Lemma \ref{le-final-add} (ii) and its proof, we know that the vertices $v_1, v_3, \dots, v_{2\alpha - 1}$ form an independent set of $G$. Thus we can assume that either $i$ and $j$ are both even, or $i$ is odd and $j$ is even.

First suppose that $d_T(v_{i+1}) = 2$. Consider $T_1 = T - v_i v_{i+1} + v_i v_j$. Clearly, $T_1$ is still a spanning tree of $G$, and $T_1$ contains one more pendent vertex than $T$, i.e., there are $n - 2\alpha + 2$ pendent vertices in $T_1$. Now from Remark \ref{remark-final-add}, we get
$$
\mu (T_1) < \alpha + \frac{1}{2}.
$$

Similarly, when $d_T(v_{j-1}) = 2$, consider $T_2 = T - v_j v_{j-1} + v_j v_i$ this time, as in last paragraph, we have
$$
\mu (T_2) < \alpha + \frac{1}{2}.
$$

Now there is still one remaining case, i.e., $d_T(v_{i+1}), d_T(v_{j-1}) \ge 3$, i.e., there are some pendent vertices attached to $v_{i+1}$ and $v_{j-1}$ in $T$. It implies that $i$ and $j$ are both even, according to the hypothesis of the structure of $T$. Consider again $T_1 = T - v_i v_{i+1} + v_i v_j$. Notice that if we delete all the pendent vertices of $T_1$, it would lead to a tree, say $T'$, of order $2\alpha-1$ obtained by identifying one end vertices of three paths $P_{k_1}$, $P_{k_2}$, $P_{k_3}$, where $k_1 + k_2 + k_3 = 2\alpha + 1$, $k_1, k_2, k_3 \ge 2$.

On one hand, if at least one of $k_1, k_2, k_3$ is even, then $T'$ has at least one pendent path of odd length, from Lemma \ref{le-final-add} (i), it infers that $T_1 \in \mathcal{T}_{t,n - t}$ with $t \le 2\alpha-2$. Now together with Lemma \ref{le-add},
we have
$$
\mu (T_1) < \alpha + \frac{1}{2}.
$$
On the other hand, if each of $k_1, k_2, k_3$ is odd, i.e., all the pendent paths of $T'$ are of even lengths, then by Lemma \ref{le-extra},
$$
\mu (T_1) < \alpha + \frac{1}{2} + \frac{1}{2(2\alpha-5)}.
$$

\noindent {\bf Case 2.} The induced subgraph of $G$ induced by the vertices $v_1, v_2, \dots, v_{2\alpha-1}$ is still $P_{2\alpha-1}$.

Note that we may assume that there is at least one pendent neighbor of $v_1$ in $T$, say $u$, which is adjacent to none of $v_2, v_3, \dots, v_{2\alpha-1}$ in $G$. Otherwise,
every neighbor of $v_1$ in $T$ is adjacent to at least one vertex of $v_2, v_3, \dots, v_{2\alpha-1}$,
we can always construct another spanning tree $T_3$, which is obtained from $P_{2\alpha-1}- v_1 = v_2 v_3 \dots v_{2\alpha-1}$ (of order $2\alpha-2$) by attaching $n - 2\alpha + 2$ pendent vertices to some vertices of $ v_2, v_3, \dots, v_{2\alpha-1}$, i.e., $T_3 \in \mathcal{T}_{2\alpha-2,n - 2\alpha + 2}$. Then by Lemma \ref{le-add},
$$
\mu (T_3) < \alpha + \frac{1}{2}.
$$

Similarly, we may also assume that there is at least one pendent neighbor of $v_{2\alpha-1}$ in $T$, say $w$, which is adjacent to none of $v_1, v_2, \dots, v_{2\alpha-2}$ in $G$.

If $u$ and $w$ are not adjacent in $G$, then $\{u, w, v_2, v_4, \dots, v_{2\alpha-2}\}$ would form an independent set of $G$, however, it contains $\alpha + 1$ vertices, which is a contradiction. Thus $u$ and $w$ are adjacent in $G$.

Denote by $a_k$ the number of pendent vertices attached to $v_{2k-1}$ in $T$, where $1 \le k \le \alpha$. Clearly,
$$
a_1 + a_2 + \dots + a_{\alpha} = n - 2 \alpha + 1.
$$
Note that we can always find some $t$, where $2 \le t \le \alpha$, such that $a_t \le \frac{n - 2 \alpha + 1}{\alpha}$, just need to assume that
$$
a_t = \min \{ a_1, a_2, \dots, a_{\alpha} \}.
$$
Here we exclude the case $t = 1$ (which is equivalent to $t = \alpha$) for convenience.
Let $T_4 = T - v_{2t - 2} v_{2t - 1} + u w$. Clearly, $T_4$ is still a spanning tree of $G$. More precisely, $T_4$ is a tree obtained from the path
$$
v_{2t-1} v_{2t} \dots v_{2\alpha - 1} w u v_1 v_2 \dots v_{2t - 3}
$$
on $2\alpha$ vertices by attaching some pendent vertices to this path. In particular, no pendent vertex outside such path on $2\alpha$ vertices is attached to $v_{2t}$ when $2 \le t \le \alpha -1$ ($w$ when $t = \alpha$, respectively) in $T_4$ (from the hypothesis of the structure of $T$). Now by Lemma \ref{le-final-2},
$$
\mu (T_4) < \alpha + \frac{1}{2} + \frac{4(\alpha-1)}{\alpha^2}
$$
follows.

Note that
$$
\frac{1}{2(2\alpha-5)} < \frac{4(\alpha-1)}{\alpha^2}.
$$
Combining Cases $1$ and $2$, by setting $T^* = T, T_1, T_2, T_3$ or $T_4$ accordingly, we would get the desired result.
\end{Proof}

By integrating the results established above, we will further improve the upper bound in Theorem \ref{old} to $\alpha + \frac{1}{2} + \frac{4(\alpha-1)}{\alpha^2}$.

\begin{Theorem} \label{th-1}
Let $G$ be a connected graph of order $n$ and independence number $\alpha$, where $n \ge 2 \alpha$ and $\alpha \ge 1$. Then $G$ has a spanning tree $T$ with
$$
\renewcommand{\arraystretch}{1.2}
\mu(T) < \left\{\begin{array}{ll}
\alpha + 1  & \mbox{if $1 \le \alpha \le 6$},\\
 \alpha + \frac{1}{2} + \frac{4(\alpha-1)}{\alpha^2}
 & \mbox{if $\alpha \ge 7$}.
\end{array}
\right.
$$
\end{Theorem}

\begin{Proof}
In Lemma \ref{le-final-add}, we have known that there exists a subtree $\bar{T}$ of $G$ of order $t$ such that the remaining $n-t$ vertices are adjacent to some vertices in $\bar{T}$, where $t \le 2 \alpha - 1$. That is to say, we can find a spanning tree $T$ of $G$ such that $T \in \mathcal{T}_{t, n-t}$, for $t \le 2 \alpha - 1$. In particular,  when $t = 2 \alpha - 1$ (i.e., $T \in \mathcal{T}_{2 \alpha - 1, n-2 \alpha + 1}$), if $\bar{T}$ is not a path, then every pendent path or internal path of $\bar{T}$ is of even length, and if $\bar{T}$ is a path, say $\bar{T} = v_1 v_2 \dots v_{2 \alpha - 1}$, then all the pendent vertices of $T$ can only be attached to the vertices among $v_1, v_3, \dots, v_{2 \alpha - 1}$.

If $t \le 2 \alpha - 2$, then by Lemma \ref{le-add}, we have
$$
\mu(T) < \alpha + \frac{1}{2}.
$$

Suppose in the following that $t = 2 \alpha - 1$. If $\bar{T}$ is not a path, then by Lemma \ref{le-extra},
$$
\mu (T) < \alpha + \frac{1}{2} + \frac{1}{2(2\alpha-5)} < \alpha + \frac{1}{2} + \frac{4(\alpha-1)}{\alpha^2}
$$
with $\alpha \ge 4$.
Next suppose that $\bar{T}$ is a path, say $\bar{T} \cong v_1 v_2 \dots v_{2\alpha - 1}$. Recall that all the pendent vertices of $T$ can only be attached to the vertices among $v_1, v_3, \dots, v_{2 \alpha - 1}$. Further, by Lemma \ref{le-key}, we can find a spanning tree of $G$, say $T^*$, such that
$$
\mu(T^*) < \alpha + \frac{1}{2} + \frac{4(\alpha-1)}{\alpha^2}.
$$
In particular,
$$
\alpha + \frac{1}{2} + \frac{4(\alpha-1)}{\alpha^2} < \alpha + 1
$$
for $\alpha \ge 7$, in which $\alpha + 1$ is the known upper bound already established in Theorem \ref{th-new}.

The proof is completed.
\end{Proof}

Similar to Corollary \ref{con}, we get the following corollary immediately.

\begin{Corollary} \label{cor-1}
Let $G$ be a connected graph of order $n$ and independence number $\alpha$, where $n \ge 2 \alpha$ and $\alpha \ge 1$. Then
$$
\renewcommand{\arraystretch}{1.2}
\mu(T) < \left\{\begin{array}{ll}
\alpha + 1  & \mbox{if $1 \le \alpha \le 6$},\\
 \alpha + \frac{1}{2} + \frac{4(\alpha-1)}{\alpha^2}
 & \mbox{if $\alpha \ge 7$}.
\end{array}
\right.
$$
\end{Corollary}

\section{Concluding remarks}

In this paper, we present an upper bound for the minimum average distance among all the spanning trees of a given connected graph, in terms of the independence number of that graph. More precisely, for a connected graph $G$ of order $n$ and independence number $\alpha$, where $n \ge 2 \alpha$, we find that there exists a spanning tree $T$ of $G$ such that
$$
\renewcommand{\arraystretch}{1.2}
\mu(T) < \left\{\begin{array}{ll}
\alpha + 1  & \mbox{if $1 \le \alpha \le 6$},\\
 \alpha + \frac{1}{2} + \frac{4(\alpha-1)}{\alpha^2}
 & \mbox{if $\alpha \ge 7$}.
\end{array}
\right.
$$
From this result, we can see that when $\alpha$ becomes larger, the upper bound on $\mu(T)$ would be better.
But it seems that we are not able to further improve the upper bound, even if $\alpha$ is rather large.

For example, let $G$ be the connected graph of order $n$ obtained from the path $P_{2k-2} = v_1 v_2 \dots v_{2k-2}$ by identifying one vertex of the complete graph on $\lfloor\frac{n-2k}{2} \rfloor + 2$ vertices with $v_1$, and identifying one vertex of the complete graph on $\lceil\frac{n-2k}{2} \rceil + 2$ vertices with $v_{2k-2}$, where $k \ge 2$ and $n \ge 2k$. It is easily seen that $H_{2k-2,n-2k+2}$ is the spanning tree of $G$ with minimum average distance. Clearly, $H_{2k-2,n-2k+2} \in \mathcal{T}_{2k - 2, n - 2k + 2}$. As we have shown in the proof of Lemma \ref{le-add},
$$
\mu(H_{2k-2,n-2k+2}) - k =  \frac{-8k ^3 + 36 k^2 - 2( 9 n + 14 ) k  + 3 n^2 + 12n }{6 n (n-1)},
$$
of which the right-side is asymptotically approximate to $\frac{1}{2}$ as $n \rightarrow \infty$, when $k$ is fixed. Note that $\alpha(G) = k$. That is to say, there exists a spanning tree $T$ of $G$ such that
$$
\mu(T) = \alpha(G) + \frac{1}{2} + o(1)
$$
as $n \rightarrow \infty$. This evidence suggests
that our new upper bound is best possible in the sense of asymptotics (when $n$ and $\alpha$ are large enough).

\section*{Statements and Declarations}

\subsection*{Funding}

This work was supported by Guangdong Basic and Applied Basic Research Foundation (Grant No.2023A1515011472), Hebei Fund for Introducing Overseas Returnees (Grant No.C20230357), Hebei Natural Science Foundation (Grant No.A2023208006), and Foundation of China Scholarship Council (Grant No.202508130088).

\subsection*{Data Availability}

Data sharing was not applicable to this article as no datasets were generated or analyzed
during the current study.

\subsection*{Competing Interests}

The authors have no relevant financial or non-financial interests to disclose.

\end{document}